\documentclass[11pt]{article} 
\usepackage{amssymb} 
\usepackage{amsmath} 
\usepackage{amsfonts} 
\usepackage{graphicx,psfrag} 
\usepackage{color} 
 
\newcommand\qed{\hfill $\square$}

\newcommand{\RR}{\mathbb R}

\newtheorem{thm}{Theorem} 
\newtheorem{prop}[thm]{Proposition} 
\newtheorem{lem}[thm]{Lemma} 
\newtheorem{cor}[thm]{Corollary}

\newcommand{\beqn}{\begin{equation}} 
\newcommand{\eeqn}{\end{equation}} 
\newcommand{\bear}{\begin{eqnarray}} 
\newcommand{\eear}{\end{eqnarray}} 
\newcommand{\bean}{\begin{eqnarray*}} 
\newcommand{\eean}{\end{eqnarray*}}

\begin{document}

\begin{center} 
{\huge \bf Localized non-diffusive  asymptotic  \\[4pt] 
patterns for  nonlinear parabolic \\[6pt]   equations with 
gradient absorption} 
\end{center}

\vspace{0.3cm}

\centerline{\large Philippe Lauren\c cot\footnote{Institut de 
Math\'ematiques de Toulouse, CNRS UMR~5219, Universit\'e Paul Sabatier 
(Toulouse~III), 118 route de Narbonne, F--31062 Toulouse Cedex 9, 
France. {\tt E-mail: laurenco@mip.ups-tlse.fr}} and Juan Luis 
V\'azquez\footnote{Departamento de Matem\'aticas, Universidad 
Aut\'onoma de Madrid, Campus de Cantoblanco, E--28049 Madrid, 
Spain. {\tt E-mail: juanluis.vazquez@uam.es}}}

\vspace{0.7cm} 
 
\rightline{\it Dedicated to Pavol Brunovsk\'y} 
 
\ 
 
 \centerline{\bf Abstract}

\medskip

{\small \noindent We study the large-time behaviour of the 
solutions $u$ of the evolution equation involving nonlinear 
diffusion and gradient absorption 
$$ 
\partial_t u - \Delta_p u + \vert\nabla u\vert^q=0. 
$$ 
We consider the problem posed for $x\in \RR^N $ and $t>0$ with 
non-negative and compactly supported initial data. We take the 
exponent $p>2$ which corresponds to slow $p$-Laplacian diffusion, 
and the exponent $q$ in the superlinear range $1<q<p-1$. In this 
range the influence of the Hamilton-Jacobi term $ \vert\nabla 
u\vert^q$ is determinant, and gives rise to the phenomenon of 
localization. The large time behaviour is described in terms of a 
suitable self-similar solution that solves a Hamilton-Jacobi 
equation. The shape of the corresponding spatial pattern is rather 
conical instead of bell-shaped or parabolic. }


\section{Introduction}\label{int}

\setcounter{thm}{0} \setcounter{equation}{0}

Researchers have been interested for decades in the long time 
description  of evolution processes where diffusion is combined 
with other effects, notably reaction, absorption and/or 
convection. The equations under study are  evolution equations of 
parabolic type, mostly nonlinear and possibly degenerate 
parabolic. A general form of such equations is 
\begin{equation}\label{geneq} 
\partial_t u =\sum_{i=1}^N \partial_{x_i} A_i(t,x,u,\nabla u) + B(t,x,u,\nabla u), 
\end{equation} 
where the $A_i$ and $B$ are nonlinear functions with suitable 
structure conditions as in \cite{LSU}, and $\nabla u$ denotes the 
spatial gradient of $u$. Under such a generality one cannot be 
very specific about the concrete asymptotic behaviour of the 
solutions as $t\to\infty$. On the contrary, when one concentrates 
into specific forms for $A_i$ and $B$, and also imposes specific 
boundary conditions, or works in the whole space for concrete 
classes of initial data, then much is known about the long time 
behaviour of the solutions, and the different types of possible 
behaviours have been carefully classified, or are in the process 
of being classified.

Following this spirit, we will concentrate here on equations where 
the right-hand side of the equation (\ref{geneq}) combines 
nonlinear diffusion of the $p$-Laplacian type with nonlinear 
absorption of Hamilton-Jacobi type. We pose the evolution on 
$\RR^N$, $N\ge 1$, as  spatial domain, and take as initial 
condition a non-negative and integrable function $u_0$ decaying 
sufficiently rapidly to zero as $|x|\to\infty$ (we may assume 
$u_0$ to be compactly supported and bounded to make the analysis 
simpler). Under these assumptions, we want to classify the 
different types of long type behaviours that may arise from the 
interaction between the two physical effects at work.  As a new 
contribution to the topic,  we describe in detail a situation in 
which the evolution leads to localized patterns with a precise 
power-like decay in time and conical spatial shape.

\subsection{The case of linear diffusion}

Before we enter into that study and in order to motivate the 
issue, it will be useful to briefly recall the main results in the 
case of standard linear diffusion. The large time behaviour of 
non-negative and integrable solutions to the diffusive 
Hamilton-Jacobi equation  (also called  viscous Hamilton-Jacobi 
equation) 
\beqn \label{a0} 
\partial_t u - \Delta u + \vert\nabla u\vert^q = 0 \;\;\mbox{ in }\;\; 
Q:=(0,\infty)\times\RR^N, 
\eeqn 
has been studied recently by 
several authors. It turns out that such behaviour is not unique 
and strongly depends on the value of the parameter $q\in 
(0,\infty)$. More precisely, always for initial data decaying 
sufficiently rapidly to zero as $|x|\to\infty$, the analysis 
reveals the existence of two critical exponents $q_1:=1$ and 
$q_*:=(N+2)/(N+1)$, and three different types of large time 
behaviour  corresponding to the parameter intervals $q\in 
(0,q_1)$, $q\in (q_1,q_*)$, and $q\in (q_*,\infty)$, respectively. 
To this we must add the study of the critical cases $q=1$ and 
$q=q_*$.

To be more precise, in the upper range $q>q_*$, the large time 
dynamics is governed by the sole diffusion term and $u$ behaves as 
a multiple of the fundamental solution to the linear heat equation 
\cite{BKL04,BGK03}. After convenient renormalization we thus see 
the typical bell-shaped Gaussian profile, in the sense that 
$$ 
t^{N/2}\ u(t,x)\to F(y)=c\ e^{-y^2/2} 
$$ 
with $y=x/t^{1/2}$, just as the same form as in the case with no 
absorption. This is a case of what is called in \cite{V91} {\sl 
asymptotic simplification}, and it also happens for the semilinear 
equation $\partial_t u - \Delta u + u^r = 0$ in $Q$ when $r> 
r_*:=(N+2)/N$.

In the intermediate range $q\in (1,q_*)$, the large time behaviour 
results from the combined effects of the diffusion and absorption 
terms and is fully described by the self-similar \textit{Very 
Singular Solution} to equation (\ref{a0}), cf.~\cite{BKL04}. The 
pattern that we see after the corresponding renormalization, 
$F_q(y)$, is a modification of the Gaussian profile corresponding 
to the Very Singular Solution and changes with $q$.

As usual in dynamical studies, the large time behaviour for the 
critical exponents $q_1$ and $q_*$ is peculiar. For $q=q_*$ it is 
still given by a multiple of the fundamental solution to the 
linear heat equation but with extra logarithmic factors resulting 
from the fact that the dynamics already feels the effects of the 
absorption term \cite{GL07}. This is usually referred to as 
\textsl{resonance}.

On the other hand, when $q$ equals $q_1=1$ then (\ref{a0}) is still a 
nonlinear equation but with the same homogeneity as a linear 
equation, and the large time behaviour does not seem to be 
thoroughly understood. Nevertheless, there are several possible 
temporal decay rates for $u$ depending on the initial condition 
\cite{BRV96,BRV97}. This is an important issue that  needs further 
investigation.

Finally, when the parameter $q$ is below $q_1$ the nonlinear 
absorption term becomes dominant and diffusion plays a secondary 
role for large times.  The typical phenomenon of the lower range 
is finite time extinction, which takes place for $q\in (0,q_1)$ 
\cite{BLSS02,Gi05}. Summarizing, the nonlinear absorption term 
only prevails for large times if $q\in (0,q_1)$. But the exponent 
$q_1=1$ is actually somehow ``doubly'' critical: indeed, not only 
does the nonlinearity rule the dynamics but it is also no longer 
locally Lipschitz continuous for $q\in (0,q_1)$, thus giving rise 
to singular phenomena such as finite time extinction.

\subsection{Nonlinear diffusion. Localized non-diffusive patterns}

In an attempt to elucidate the true role of the absorption term 
when it governs the dynamics, we devote this paper to investigate 
in detail the large time behaviour of non-negative solutions $u$ 
to the related Cauchy problem 
\bear \label{a1} 
\partial_t u - \Delta_p u + \vert\nabla u\vert^q & = & 0\ \,, \quad 
(t,x)\in Q\,,\\ 
\label{a2} u(0) & = & u_0\ge 0\,, \quad x\in\RR^N\,, 
\eear 
where the linear diffusion operator $\Delta u$ is replaced by the 
$p$-Laplacian operator 
\beqn \label{a3} \Delta_p u := \mbox{ div } 
\left( \vert\nabla u\vert^{p-2}\ \nabla u \right)\,. 
\eeqn 
When $p>2$, (\ref{a1}) is a quasilinear degenerate parabolic equation 
which reduces to the semilinear equation (\ref{a0}) for $p=2$, and has 
been studied recently: the existence and properties of 
self-similar solutions are investigated in \cite{GB02,Sh04a}, 
while qualitative properties of non-negative and integrable 
solutions are obtained in \cite{ATU04,BtLxx,Y96}. In 
particular, two critical exponents have been recently identified 
for (\ref{a1}) \cite{ATU04,BtLxx}, namely 
$$ 
q_1:= p-1 \;\;\mbox{ and }\;\; q_* := p - \frac{N}{N+1}\,. 
$$ 
As for $p=2$, the nonlinear absorption term governs the large time 
dynamics for $q\in (0,q_1)$ but we now have $q_1>1$ since $p>2$. 
However, the singular phenomena like extinction in finite 
time are only expected to happen in a still lower subregion $q\in (0,1)$.

Consequently, in the new intermediate range  $q\in (1,q_1)$ that 
we want to explore here the nonlinearity is locally Lipschitz 
continuous and no singular phenomenon (such as extinction) can 
occur. Still, the absorption term strongly affects the evolution 
for $q\in (1,q_1)$. Two results  are established in \cite{BtLxx} 
in that parameter interval  for compactly supported initial data: 
there are positive constants $C>0$ and $R>0$ (possibly depending 
on the initial data) such that 
$$ 
{\cal P}_u(t) \subset B(0,R) \;\;\mbox{ and }\;\; \Vert 
u(t)\Vert_\infty \le C\ t^{-1/(q-1)} 
$$ 
for all $t\ge 0$ [here and below we use the notation ${\cal P}_u(t)=\left\{ 
x\in\RR^N:  \;\; u(t,x)>0 \right\} $ for the positivity set of the 
function $u$ at time $t$ and we use the abbreviated notation 
$u(t)$ instead of $u(t,\cdot)$ and so on, whenever there is no fear of confusion].

These properties are reminiscent of those enjoyed by the viscosity 
solutions to the Hamilton-Jacobi equation 
\beqn \label{hje} 
\partial_t h + \vert\nabla h\vert^q = 0 \;\;\mbox{ in }\;\; Q 
\eeqn 
with compactly supported initial data. 
Indeed, denoting by $h$ the viscosity solution to equation 
(\ref{hje}) with initial condition $h(0)=u(0)$, we have for each 
$t\ge 0$ (cf.  \cite{Bl85}) 
$$ 
{\sup_{t\ge 0}{\left\{ t^{1/(q-1)}\  \Vert h(t)\Vert_\infty 
\right\}} < \infty}\,, \quad \mbox{and } \quad {\cal P}_h(t)= 
{\cal P}_h(0). 
$$ 
This is in sharp contrast with the behaviour of the solution $w$ 
to the  pure  $p$-Laplacian equation $\partial_t w - \Delta_p w = 
0$ in $(0,\infty)\times\RR^N$ with the same initial condition 
$w(0)=u(0)$. Indeed, the $L^\infty$-norm of $w(t)$ only decays at 
the much slower rate $$ \|w(t)\|_\infty=O(t^{-N/(N(p-2)+p)}), 
$$ and the 
positivity set ${\cal P}_w(t)=\left\{ x\in\RR^N : \; w(t,x)>0 
\right\}$ of $w(t)$ expands and  as $t\to\infty $ it fills the 
whole space  $\RR^N$. The absorption term $|\nabla u|^q$ in 
(\ref{a1}) thus prevents this expansion if $q\in (1,q_1)$ and the 
positivity set of $u(t)$ remains \textit{localized} in a fixed 
ball of $\RR^N$ for all times.

The property of localization is already known to be satisfied by 
compactly supported non-negative solutions to second-order 
degenerate parabolic equations with a sufficiently strong 
absorption involving the solution only as, for instance, 
$\partial_t w - \Delta_p w + w^r = 0$ in $Q$ 
when $r\in (1,p-1)$,  cf. \cite{DV85,Ka87,Y96}. It had apparently 
remained unnoticed for second-order degenerate parabolic equations 
with an absorption term depending solely on the gradient.

According to the previous discussion, the qualitative results 
obtained so far on compactly supported non-negative solutions to 
(\ref{a1}) show evidence of the domination of the nonlinear 
absorption term for large times. The purpose of this paper is to 
go one step further in that direction and show that $u(t)$ behaves 
like a self-similar solution to (\ref{hje}) as $t\to\infty$.

For comparison, it is interesting to notice that this localized 
asymptotics is quite different in shape from the behaviour of the 
diffusion-absorption equation \ $\partial_t w - \Delta w^m + w^r = 
0$ in $Q$ that offers a strong similarity in 
other respects with our present problem when $1<r<m$ . The 
corresponding localized behaviour has been studied in \cite{BNP, BKP, ChV99} 
and leads to a mesa-like pattern. In this 
paper we show that a conical pattern $V_\infty$ is gradually 
formed, precisely given by  formula \eqref{a10}, 
and is actually a viscosity solution to the stationary equation 
\beqn 
\vert\nabla v\vert^q - v = 0 \;\;\mbox{ in }\;\; \RR^N\,. 
\eeqn 
It is associated to an exact self-similar viscosity solution 
of the Hamilton-Jacobi equation (\ref{hje}).

As a final contribution, we  investigate conditions under which the support 
does not move at all, so-called {\sl infinite waiting time}, so 
that what happens for the whole time span is a process of internal 
reorganization. We refer to \cite{Ar85} and \cite{V07} for waiting 
times in porous medium flows. In those cases the waiting time is 
always finite. An infinite waiting time is described in 
\cite{ChV96, ChV99} for the equation $\partial_t w - \Delta w^m + 
w^r = 0$ in $Q$ when $1<r<m$.  Properties of localization and waiting 
times for more general equations of the form \eqref{geneq} are 
considered for instance in \cite{Ant96, AntDiaz89, DV85}.


\section{Preliminaries and main results}

Before stating our results, let us first specify our assumptions 
and recall the properties of solutions to (\ref{a1}), (\ref{a2}) 
with non-negative and compactly supported initial data established 
in \cite{BtLxx}. We assume that \beqn \label{a4} p>2 \;\;\mbox{ 
and }\;\; q\in (1,p-1)\,, \eeqn and that the initial condition 
$u_0$ enjoys the following properties: 
\beqn \label{a5} 
u_0\in W^{1,\infty}(\RR^N)\,, \;\; u_0\ge 0\,, \;\; 
\mathcal{P}_0:=\left\{ x\in\RR^N : \;\; u_0(x)>0 \right\} \subset B(0,R_0)\,, 
\;\; u_0\not\equiv 0\,, \eeqn 
for some $R_0>0$. Let us also introduce an exponent that will play a role 
in what follows: $\xi:=1/(q(N+1)-N)>0$. It goes from $\xi=1$ for 
$q=1$ to $\xi=1/(N(p-2)+p-1)$ when $q=p-1$. 
\begin{prop}\label{pra1} 
Under the above assumptions, the Cauchy problem {\rm (\ref{a1}), 
(\ref{a2})} has a unique non-negative viscosity solution 
$$ 
u\in\mathcal{BC}([0,\infty)\times\RR^N)\cap 
L^\infty(0,\infty;W^{1,\infty}(\RR^N)) $$ which satisfies: 
\begin{itemize} 
\item[(i)] There is $R_\infty>0$ such that \beqn \label{a6} \mbox{ 
supp }(u(t)) \subset B(0,R_\infty) \;\;\mbox{ for all }\;\; t\ge 
0\,. \eeqn \item[(ii)] There are positive constants $C_0$ and 
$C_1$ such that \beqn \label{a7} \left\{ 
\begin{array}{l} 
\Vert u(t)\Vert_\infty \le C_0\ \Vert u(s)\Vert_1^{q\xi}\ 
(t-s)^{-N\xi}\,, \\ \\ 
\Vert\nabla u(t)\Vert_\infty \le C_0\ \Vert u(s)\Vert_1^{\xi}\ 
(t-s)^{-(N+1)\xi}\,, \end{array} \right. \eeqn \beqn \label{a8a} 
\Vert u(t)\Vert_1 \le \Vert u_0\Vert_1\,, \quad \Vert 
u(t)\Vert_\infty \le \Vert u_0\Vert_\infty\,, \quad  \Vert \nabla 
u(t)\Vert_\infty \le \Vert \nabla u_0\Vert_\infty\,, \eeqn \beqn 
\label{a8b} t^{1/(q-1)}\ \left( \Vert u(t)\Vert_1 + \Vert 
u(t)\Vert_\infty + \Vert \nabla u(t)\Vert_\infty \right) \le C_1 
\eeqn for all $t>0$ and $s\in [0,T)$. 
\end{itemize} 
\end{prop}

Here and below $\mathcal{BC}([0,\infty)\times\RR^N)$ denotes the 
space of bounded and continuous functions on 
$[0,\infty)\times\RR^N$ and $\Vert \cdot\Vert_r$ denotes the 
$L^r(\RR^N)$-norm for $r\in [1,\infty]$.

Assertions (\ref{a6}), (\ref{a7}), and (\ref{a8a}) are proved in 
\cite[Theorem~1.5~(i)]{BtLxx} and \cite[Proposition~1.3]{BtLxx}, 
respectively, while (\ref{a8b}) is established in 
\cite[Corollary~1.6~(i)]{BtLxx} for the $L^1$-norm. The estimates 
(\ref{a8b}) for $\Vert u(t)\Vert_\infty$ and $\Vert \nabla 
u(t)\Vert_\infty$  follow then from that for $\Vert u(t)\Vert_1$ 
by (\ref{a7}) with $s=t/2$.

\subsection{Statement of the main asymptotic results}

As indicated above we call $\mathcal{P}_u(t)$ the 
\textit{positivity set} of the solution $u$ at time $t\ge 0$, that 
we abbreviate as $ \mathcal{P}(t) $. Since $u(t)$ is compactly 
supported and continuous for each $t\ge 0$ by 
Proposition~\ref{pra1}, $\mathcal{P}(t)$ is a  bounded open subset 
of $\RR^N$ for all $t\ge 0$. Much more can actually be said on the 
family $\{ \mathcal{P}(t) \}_{t\ge 0}$ and it will  allow us to 
identify the large time behaviour of $u$. 
\begin{thm}\label{tha2} 
The mapping $t\longmapsto \mathcal{P}(t)$ gives a non-decreasing 
family of bounded open subsets of $\RR^N$ and moreover 
$$ 
\mathcal{P}_\infty := \bigcup_{t\ge 0} \mathcal{P}(t) \;\;\mbox{ 
is a bounded open subset of }\;\; \RR^N\,. 
$$ 
In addition, 
\beqn 
\label{a9} \lim_{t\to\infty} \left\Vert t^{1/(q-1)}\ u(t) - 
V_\infty \right\Vert_\infty = 0\,, \eeqn 
where 
\beqn \label{a10} 
V_\infty(x) := \frac{q-1}{q^{q/(q-1)}}\ \left( 
\inf_{y\in\RR^N\setminus\mathcal{P}_\infty}{\{ \vert y-x\vert\}} 
\right)^{q/(q-1)} \;\;\mbox{ for }\;\; x\in\RR^N\,. \eeqn 
\end{thm}

\begin{figure}\label{figolu} 
\centerline{\includegraphics[height=6cm]{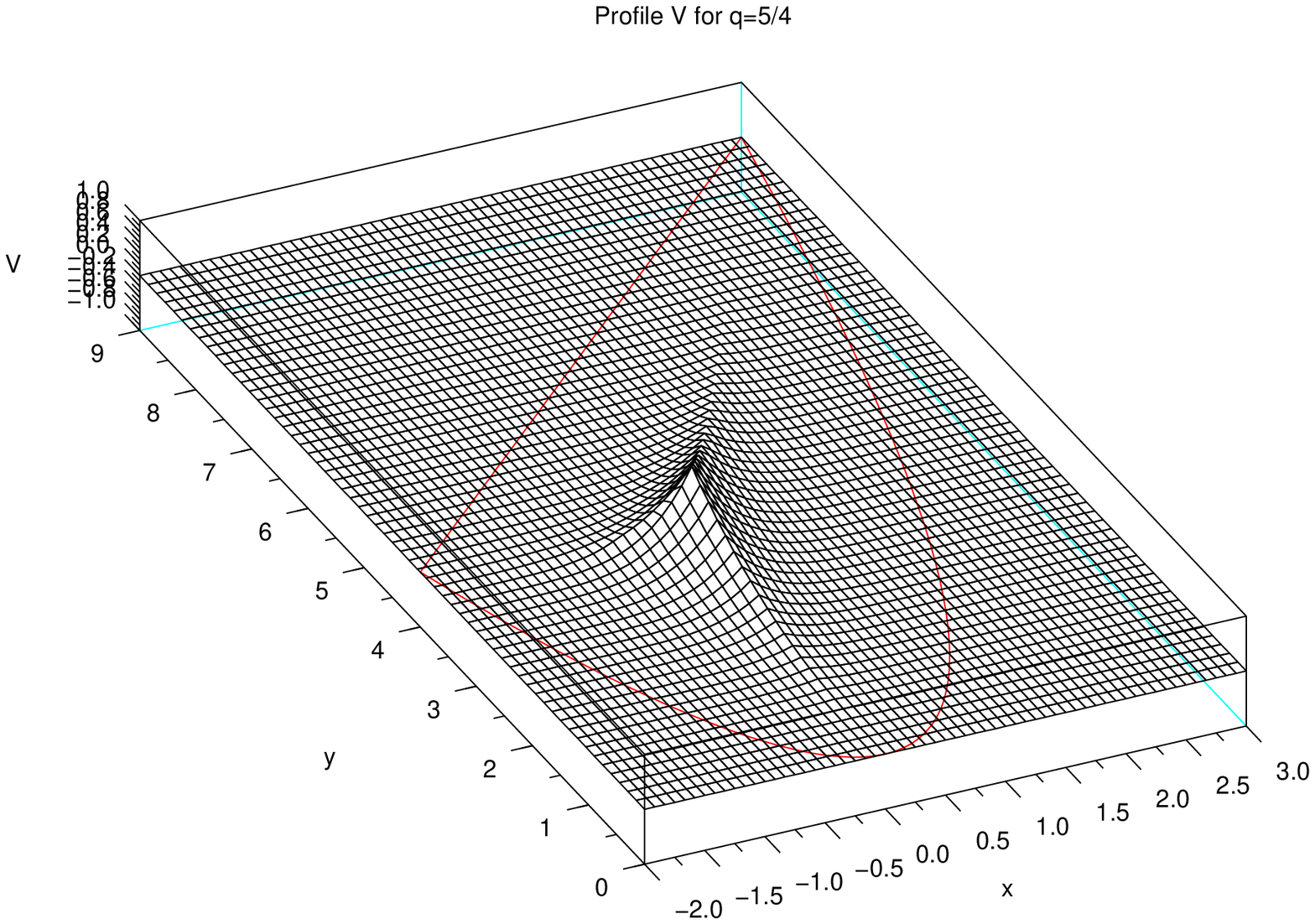} \hskip0.3cm 
            \includegraphics[height=6cm]{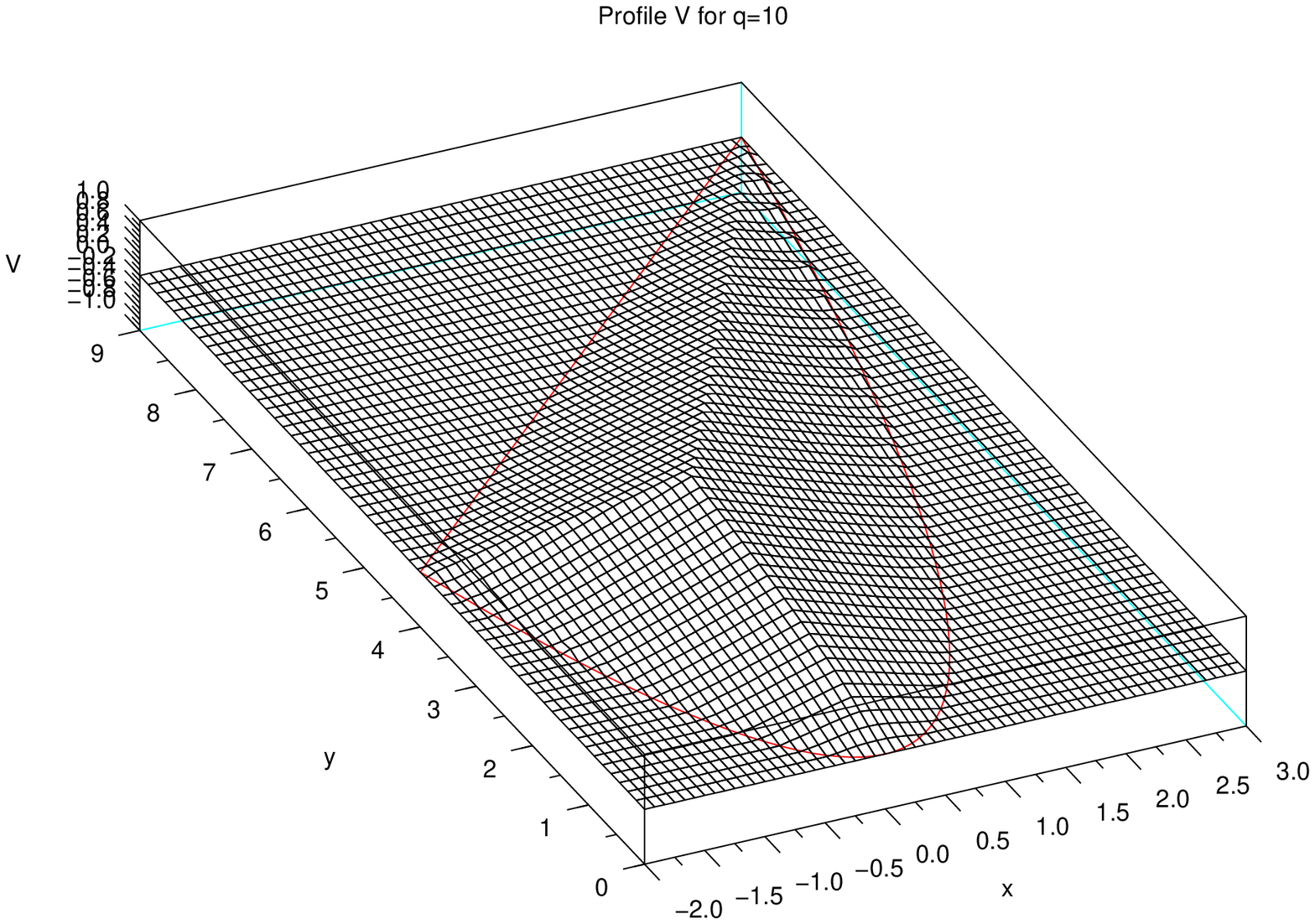}} 
\caption{Profile $V_\infty$ for $q=5/4$ (left) and $q=10$ (right)} 
\end{figure}

We first point out that $U(t,x):= t^{-1/(q-1)}\ V_\infty(x)$ is a 
self-similar (viscosity) solution to (\ref{hje}) which vanishes 
outside $\mathcal{P}_\infty$ and is positive in 
$\mathcal{P}_\infty$. Therefore, (\ref{a9}) asserts that the 
diffusion term no longer matters for large times as claimed 
previously. Still, it might influence the shape of 
$\mathcal{P}_\infty$ during the time evolution. We will come back 
to this issue in Theorems~\ref{tha3} and~\ref{tha4} below. Also, 
the final (rescaled) pattern $V_\infty$ has a conical shape as 
illustrated in Figure~$1$ for $N=2$ and $\mathcal{P}_\infty=\{ 
(x,y)\in [-2,3]\times [0,9] : x^2 < y < x+6 \}$.

The proof of Theorem~\ref{tha2} is divided in several steps and is 
performed in the forthcoming Sections~\ref{tmps} and~\ref{css}. 
We nevertheless outline its main steps now: we 
first establish the monotonicity property of the family $\{ 
\mathcal{P}(t) \}_{t\ge 0}$ by comparison arguments, and lower 
bounds for $u(t,x)$ as well. We next turn to the proof of 
(\ref{a9}) and introduce the \textit{self-similar} variables; we 
keep the space variable $x$ (since we expect localized behaviour), 
and (as usual) introduce logarithmic time 
$$ 
\tau:=\ln{(1+(q-1)t)}/(q-1), 
$$ 
as well as  the new unknown function $v=v(\tau,x)$ defined by 
\beqn \label{a11} 
u(t,x) = 
(1+(q-1)t)^{-1/(q-1)}\ v\left( \tau, x \right)\,, \quad (t,x)\in 
[0,\infty)\times\RR^N\,. 
\eeqn 
Usually, $v(\tau)=v(\tau,\cdot)$ is called the 
renormalized trajectory. Clearly, $v$ solves 
\bear \label{a12} 
\partial_\tau v + \vert\nabla v\vert^q -v & = & e^{-(p-1-q)\tau}\ 
\Delta_p v\ \,, 
\quad (\tau,x)\in Q\,,\\ 
\label{a13} v(0) & = & u_0\,, \quad x\in\RR^N\,. \eear The 
convergence (\ref{a9}) then reads $\|v(\tau) - 
V_\infty\|_\infty\longrightarrow 0$ as $\tau\to\infty$ and we are 
thus left to establish the convergence of $v(\tau)$ as 
$\tau\to\infty$. We first notice that, as $q<p-1$, the coefficient 
in front of the diffusion term vanishes for large times, so that 
the diffusion term is clearly negligible for large times. To 
identify the large time behaviour of $v(\tau)$ we use the 
half-relaxed limits technique \cite{BlP88} in the same fashion as 
in \cite[Theorem~1]{NR99} and \cite[Section~3]{Ro01}, and 
introduce 
$$ 
v_*(\tau,x):= 
\liminf_{(\sigma,y,s)\to(\tau,x,\infty)}{v(\sigma+s,y)} \;\;\mbox{ 
and }\;\; v^*(\tau,x):= 
\limsup_{(\sigma,y,s)\to(\tau,x,\infty)}{v(\sigma+s,y)} $$ for 
$(\tau,x)\in Q$. We first prove that $v_*$ and 
$v^*$ do not depend on $\tau$ and are viscosity supersolution and 
subsolution, respectively, to the Hamilton-Jacobi equation 
$|\nabla z|^q-z=0$ in $\RR^N$. Noticing that $v_*$ and $v^*$ are 
positive on the same open subset $\mathcal{P}_\infty$ of $\RR^N$, 
a comparison argument can be applied to conclude that $v_*=v^*$, 
from which the convergence of $v(\tau)$ follows. 
\medskip

\subsection{The question of waiting times}

In Theorem~\ref{tha2} we have identified the first term of the 
expansion of $u(t)$ as $t\to\infty$. Still there is an unspecified 
quantity, namely the set $\mathcal{P}_\infty$ on which additional 
information would be welcome. In particular, is there a class of 
initial data for which the positivity set does not change during 
time evolution? The next result gives a positive answer to this 
question. 
\begin{thm}\label{tha3} 
Assume that there exists $\delta>0$ such that \beqn \label{a14} 
u_0(x) \le a_0\ |x-x_0|^{(p-q)/(p-1-q)} \;\;\mbox{ for }\;\; x\in 
B(x_0,\delta) \;\;\mbox{ and }\;\; x_0\in\partial\mathcal{P}_0\,, 
\eeqn \beqn \label{a15} \| u_0\|_\infty \le a_0\ 
\delta^{(p-q)/(p-q-1)}\,, \eeqn where \beqn \label{a16} a_0 := 
\frac{p-1-q}{p-q}\ \left( N - 1 + \frac{p-1}{p-1-q} 
\right)^{-1/(p-1-q)}\,. \eeqn Then $\mathcal{P}(t)=\mathcal{P}_0$ 
for every $t\ge 0$, and thus $\mathcal{P}_\infty=\mathcal{P}_0$. 
\end{thm}

When a point of the boundary of the positivity set stays as such 
for a time we say that there is a waiting time at that point. 
Waiting times are typical in degenerate parabolic equations with 
slow diffusion, like $\partial_t w=\Delta w^m$ for $m>1$ or 
$\partial_t w=\Delta_p w$ for $p>2$, cf. \cite{V07}, but in these 
cases they are finite. In the present case we have exhibited an 
infinite waiting time at all points of the boundary of the initial 
support.

For small values of $q$, the behaviour reported in 
Theorem~\ref{tha3} ceases to be observed as soon as (\ref{a14}) is 
not fulfilled, showing the optimality of the exponent 
$(p-q)/(p-1-q)$ for the support not to evolve.

\begin{thm}\label{tha4} 
Assume that $x_0\in\partial\mathcal{P}_0$ and there is $\delta>0$ 
such that \beqn \label{a17} u_0(x) \ge A\ |x-x_0|^{(p-q)/(p-1-q)} 
\;\;\mbox{ for }\;\; x\in B(x_0,\delta)\,. \eeqn If 
\beqn 
\label{a18} 
q< q_2:=\min{\left\{ p- \frac{2N}{N+1} , \frac{p}{2} 
\right\}} 
\eeqn 
there is $A_q>0$ depending only on $N$, $p$, and 
$q$ such that, if $A\ge A_q$, 
$$ 
u(t,x_0)>0 \;\;\mbox{ for all }\;\; t>0\,. 
$$ 
In particular, $x_0\in\mathcal{P}(t)$ for each $t>0$. 
\end{thm}

A similar result is actually also valid for $q\in [q_2,p-1)$ but 
requires a stronger assumption than (\ref{a17}), see 
Proposition~\ref{prd3} below. 

\section{Time monotonicity of the positivity set}\label{tmps} 

We first establish the time monotonicity of $\mathcal{P}(t)$.

\begin{prop} 
\label{prb2} 
For $t_1\in [0,\infty)$ and $t_2\in (t_1,\infty)$ we have 
\ $\mathcal{P}(t_1) \subseteq \mathcal{P}(t_2)$ 
and 
\beqn 
\label{b2} 
\mathcal{P}_\infty := \bigcup_{t\ge 0} \mathcal{P}(t) \;\;\mbox{ is a 
bounded open subset of }\;\; \RR^N\,. 
\eeqn 
\end{prop}

The proof relies on the construction of a suitable subsolution which 
we perform next.

\begin{lem}\label{leb1} 
Fix $R>0$ and put 
$$ 
s_A(t) := A\ (1+t)^{-1/(q-1)}\ \left( R^2 - \vert x\vert^2 
\right)^{q/(q-1)}\,, \quad (t,x)\in [0,\infty)\times B(0,R)\,. 
$$ 
There exists $A_R>0$ such that $s_A$ is a (viscosity) subsolution to 
{\rm (\ref{a1})} in $(0,\infty)\times B(0,R)$ for every $A\in [0,A_R]$. 
\end{lem}

\noindent\textbf{Proof.} Introducing $\sigma(x):= \left( R^2 - |x|^2 
\right)^{q/(q-1)}$ for $x\in B(0,R)$, we observe that $\sigma$ and 
$|\nabla\sigma|^{p-2} \nabla\sigma$ both belong to 
$\mathcal{C}^1(B(0,R))$. Therefore, if $(t,x)\in [0,\infty)\times 
B(0,R)$, 
\bean 
& & \frac{(1+t)^{q/(q-1)}}{A\ \sigma(x)}\ \left\{ \partial_t s_A(t,x) 
- \Delta_p s_A(t,x) + |\nabla s_A(t,x)|^q \right\}\\ 
& \le & - \frac{1}{q-1} + \left( \frac{2q}{q-1} \right)^q\ A^{q-1}\ 
|x|^q \\ 
& & + \left( \frac{2q}{q-1} \right)^q\ \frac{(N+p-2)\ 
A^{p-2}}{(1+t)^{(p-1-q)/(q-1)}}\ \sigma(x)^{(p-1-q)/q}\ |x|^{p-2}\\ 
& \le & - \frac{1}{q-1} + \left( \frac{2qR}{q-1} \right)^q\ A^{q-1} \\ 
& & + \left( \frac{2q}{q-1} \right)^q\ (N+p-2)\ A^{p-2}\ 
R^{(2(p-1-q)/(q-1))+p-2} \le  0 \eean if $A\in [0,A_R]$ where 
$A_R$ is sufficiently small and depends only on $N$, $p$, $q$, and 
$R$. \qed

\medskip

\noindent\textbf{Proof of Proposition~\ref{prb2}.} Fix $t_1\ge 0$ and 
$x_1\in \mathcal{P}(t_1)$. Owing to the continuity of $x\mapsto 
u(t,x)$ there are $\delta>0$ and $R>0$ such that $u(t_1,x)\ge\delta>0$ 
for $x\in B(x_1,R)$. Choosing $A\in (0,A_R)$ such that $A\le\delta\ 
R^{-2q/(q-1)}$, we notice that 
$$ 
s_A(t_1,x-x_1) \le A\ R^{2q/(q-1)} \le \delta \le u(t_1,x) \;\;\mbox{ 
for }\;\; x\in B(x_1,R)\,, 
$$ 
and 
$$ 
s_A(t,x-x_1) = 0 \le u(t,x) \;\;\mbox{ for }\;\; (t,x)\in 
[t_1,\infty)\times\partial B(x_1,R)\,, 
$$ 
the parameter $A_R$ and the function $s_A$ being defined in 
Lemma~\ref{leb1}. By Lemma~\ref{leb1}, $(t,x)\longmapsto s_A(t,x-x_1)$ 
is a subsolution to (\ref{a1}) in $[t_1,\infty)\times B(0,R)$ and we 
infer from the comparison principle \cite[Theorem~8.2]{CIL92} that 
$s_A(t,x-x_1) \le u(t,x)$ for $(t,x)\in [t_1,\infty)\times B(0,R)$. In 
particular, $u(t,x_1)\ge s_A(t,0)= A\ (1+t)^{-1/(q-1)}\ R^{2q/(q-1)} > 
0$ for $t\ge t_1$. Therefore, $x_1\in\mathcal{P}(t)$ for every $t\ge 
t_1$, which proves the first assertion of Proposition~\ref{prb2}. The 
set $\mathcal{P}_\infty$ is then clearly an open subset of $\RR^N$ 
which is bounded by (\ref{a6}). \qed

\medskip

Another useful consequence of Lemma~\ref{leb1} is a lower bound for 
$u(t)$ which implies in particular that $\|u(t)\|_\infty$ cannot decay 
to zero at a faster rate than $t^{-1/(q-1)}$ for large times.

\begin{cor} 
\label{cob3} 
For each $x\in\mathcal{P}_\infty$ there are $T_x\ge 0$ and 
$\varepsilon_x>0$ such that 
\beqn 
\label{b3} 
u(t,x) \ge \varepsilon_x\ (1+t)^{-1/(q-1)} \;\;\mbox{ for }\;\; t\ge 
T_x\,. 
\eeqn 
In particular, there exists $C_2>0$ such that 
\beqn 
\label{b3b} 
\Vert u(t)\Vert_\infty \ge C_2\ t^{-1/(q-1)} \;\;\mbox{ for }\;\; t\ge 
1\,. 
\eeqn 
\end{cor}

\noindent\textbf{Proof.} Consider $x\in\mathcal{P}_\infty$. Then there 
is $T_x\ge 0$ such that $x\in\mathcal{P}(t)$ for $t\ge T_x$. The 
continuity of $y\longmapsto u(T_x,y)$ ensures that there are 
$\delta>0$ and $R>0$ such that $u(T_x,y)\ge\delta>0$ for $y\in 
B(x,R)$. Choosing $A\in (0,A_R)$ such that $A\le \delta\ 
R^{-2q/(q-1)}$ we argue by comparison as in the proof of 
Proposition~\ref{prb2} to conclude that $s_A(t,y-x) \le u(t,y)$ for 
$(t,y)\in [T_x,\infty)\times B(x,R)$ (recall that $A_R$ and $s_A$ are 
defined in Lemma~\ref{leb1}). Consequently, $u(t,x)\ge s_A(t,0)= A\ 
(1+t)^{-1/(q-1)}\ R^{2q/(q-1)} > 0$ for $t\ge T_x$, whence (\ref{b3}) 
with $\varepsilon_x=A\ R^{2q/(q-1)}$. The lower bound (\ref{b3b}) is 
then a straightforward consequence of (\ref{b3}) with any 
$x\in\mathcal{P}_0$. \qed

\section{Convergence to self-similarity}\label{css}

We now investigate the properties and convergence for large 
times of the function $v$ defined in (\ref{a11}). For $\tau\ge 0$ 
let $\mathcal{P}_v(\tau):= \left\{ x\in\RR^N: \;\; v(\tau,x)>0 
\right\}$ be the positivity set of $v(\tau)$. According to 
(\ref{a11}), $\mathcal{P}_v(\tau)$ is given by 
$$ 
\mathcal{P}_v(\tau) = \mathcal{P}_u\left( \frac{e^{(q-1)\tau}-1}{q-1} 
\right) \;\;\mbox{ for }\;\; \tau\ge 0\,. 
$$

We next state some bounds for the rescaled function $v$ which are easy 
consequences of (\ref{a6}), (\ref{a8a}), (\ref{a8b}), and (\ref{b3}).

\begin{lem} 
\label{lec1} 
There is $C_3>0$ such that 
\beqn 
\label{c1} 
\Vert v(\tau)\Vert_\infty + \Vert \nabla v(\tau)\Vert_\infty \le C_3 
\;\;\mbox{ for }\;\; \tau\ge 0\,. 
\eeqn 
In addition, for each $x\in\mathcal{P}_\infty$, there are $\tau_x$ and 
$\varepsilon_x>0$ such that 
\beqn 
\label{c1b} 
v(\tau)\ge \varepsilon_x>0 \;\;\mbox{ for }\;\; \tau\ge \tau_x\,. 
\eeqn

Finally, $\tau\longmapsto \mathcal{P}_v(\tau)$ is a non-decreasing 
family of bounded open subsets of $\RR^N$ and 
$$ 
\bigcup_{\tau\ge 0} \mathcal{P}_v(\tau) = \mathcal{P}_\infty\,, 
$$ 
the limit set $\mathcal{P}_\infty$ being defined by {\rm (\ref{b2})}. 
\end{lem}

\noindent\textbf{Proof.} The estimates (\ref{c1}) readily follow from 
(\ref{a8a}), (\ref{a8b}), and (\ref{a11}), while (\ref{c1b}) and the 
properties of $\{\mathcal{P}_v(\tau)\}_{\tau\ge 0}$ are 
straightforward consequences of Proposition~\ref{prb2} and 
Corollary~\ref{cob3}. \qed

\medskip

We turn to our main task, i.\,e., the behaviour of $v(\tau)$ as $\tau\to\infty$, and 
actually aim at showing the convergence of $v(\tau)$ as 
$\tau\to\infty$ to an asymptotic profile. For that purpose we use the half-relaxed limits 
technique \cite{BlP88} in the same fashion as in 
\cite[Section~3]{Ro01}. For $(\tau,x)\in [0,\infty)\times\RR^N$ we 
define 
\beqn 
\label{c2} 
v_*(x):= \liminf_{(\sigma,y,s)\to(\tau,x,\infty)}{v(\sigma+s,y)} 
\;\;\mbox{ and }\;\; v^*(x):= 
\limsup_{(\sigma,y,s)\to(\tau,x,\infty)}{v(\sigma+s,y)} 
\eeqn 
and first note that the right-hand sides of the above definitions 
indeed do not depend on $\tau\ge 0$. In addition, 
\beqn 
\label{c3} 
0 \le v_*(x) \le v^*(x) \;\;\mbox{ for }\;\; x\in\RR^N 
\eeqn 
by (\ref{c2}), while Lemma~\ref{lec1} and the Rademacher theorem 
clearly ensure that $v_*$ and $v^*$ both belong to 
$W^{1,\infty}(\RR^N)$. Finally, by \cite[Th\'eor\`eme~4.1]{Bl94} 
applied to equation (\ref{a12}), $v^*$ and $v_*$ are viscosity 
subsolution and supersolution, respectively, to the Hamilton-Jacobi equation 
\beqn 
\label{c4} 
H(z,\nabla z) := |\nabla z|^q - z = 0 \;\;\;\mbox{ in }\;\;\; 
\RR^N\,. 
\eeqn

The next step is to show that $v^*$ and $v_*$ actually coincide. At 
this point we emphasize that $v^*$ and $v_*$ enjoy an additional 
property, namely, 
\beqn 
\label{c5} 
\left\{ 
\begin{array}{lcl} 
v^*(x) = v_*(x) = 0 & \mbox{ for } & 
x\in\RR^N\setminus\mathcal{P}_\infty\,, \\ 
 & & \\ 
v^*(x) \ge v_*(x)> 0 & \mbox{ for } & x\in \mathcal{P}_\infty 
\end{array} 
\right. 
\eeqn 
(a similar situation is encountered in the proof of 
\cite[Theorem~1]{NR99}). 
Indeed, if $x\not\in\mathcal{P}_\infty$, the time monotonicity of 
$\tau\mapsto \mathcal{P}_v(\tau)$ (Lemma~\ref{lec1}) warrants that 
$x\not\in\mathcal{P}_v(\tau)$ for all $\tau\ge 0$. Then, if 
$s_n\to\infty$, $\sigma_n\to 0$, and $x_n\to x$ are such that 
$v(\sigma_n+s_n,x_n)\longrightarrow v^*(x)$ as $n\to \infty$, it 
follows from (\ref{c1}) that 
$$ 
v(\sigma_n+s_n,x_n) \le v(\sigma_n+s_n,x) + |x-x_n|\ \|\nabla 
v(\sigma_n+s_n)\|_\infty \le C_3\ |x-x_n|\,. 
$$ 
Passing to the limit and recalling (\ref{c3}) gives the first 
assertion in (\ref{c5}). Consider next $x\in\mathcal{P}_\infty$ and 
recall that there are $\tau_x$ and $\varepsilon_x>0$ such that 
$v(\tau,x)\ge\varepsilon_x$ for $\tau\ge\tau_x$ by (\ref{c1b}). Pick 
sequences $(\sigma_n)_{n\ge 1}$, $(s_n)_{n\ge 1}$, and $(x_n)_{n\ge 
1}$ such that $\sigma_n\to 0$, $s_n\to\infty$, $x_n\to x$, and 
$v(\sigma_n+s_n,x_n)\longrightarrow v_*(x)$ as $n\to\infty$. For $n$ 
large enough we have $\sigma_n+s_n\ge\tau_x$ and we infer from 
Lemma~\ref{lec1} that 
$$ 
v(\sigma_n+s_n,x_n) \ge v(\sigma_n+s_n,x) - C_3\ |x-x_n| \ge 
\varepsilon_x - C_3\ |x-x_n|\,. 
$$ 
Letting $n\to\infty$ gives $v_*(x)\ge\varepsilon_x>0$ and completes 
the proof of (\ref{c5}).

We next introduce 
\beqn 
\label{c6} 
V_*(x) := \frac{q}{q-1}\ v_*(x)^{(q-1)/q} \;\;\mbox{ and }\;\; V^*(x) 
:= \frac{q}{q-1}\ v^*(x)^{(q-1)/q} 
\eeqn 
for $x\in\mathcal{P}_\infty$. Arguing as in 
\cite[Corollaire~2.1]{Bl94} or \cite[Proposition~II.2.5]{BdCD97}, it 
easily follows from the properties of $v^*$ and $v_*$ that $V^*$ and 
$V_*$ are viscosity subsolution and supersolution, respectively, to 
the eikonal equation 
\beqn 
\label{c7} 
H_{ei}(z) := |\nabla z|- 1 = 0 \;\;\mbox{ in }\;\; 
\mathcal{P}_\infty\,. 
\eeqn 
Now, $H_{ei}$ depends solely on $\nabla z$ and is a convex function 
of $\nabla z$, so that the assumptions (H1), (H2), and (H4) in 
\cite{I87} are clearly fulfilled. Furthermore, $H_{ei}(0)=-1<0$ which 
warrants that the assumption (H3) in \cite{I87} is also fulfilled 
(the function $\varphi$ in \cite{I87} being here identically zero). 
Since $V_*=V^*=0$ on 
$\partial\mathcal{P}_\infty$, we are in a position to apply 
\cite[Theorem~1]{I87} to conclude that $V^*(x)\le V_*(x)$ for $x\in 
\mathcal{P}_\infty$. Recalling (\ref{c5}) and (\ref{c6}) we end up 
with $V^*(x)=V_*(x)$ for $x\in \mathcal{P}_\infty$ and thus 
$v_*(x)=v^*(x)$ in $\RR^N$. Owing to \cite[Lemme~4.1]{Bl94} or 
\cite[Lemma~V.1.9]{BdCD97}, the equality $v_*=v^*$ and (\ref{c2}) 
provide the convergence of $(v(s))_{s\ge 1}$ towards $v_*$ uniformly 
on every compact subset of $\RR^N$ as $s\to\infty$. Since all these 
functions are compactly supported in the closure of 
$\mathcal{P}_\infty$, we conclude that 
$$ 
\lim_{\tau\to\infty} \|v(\tau)-v_*\|_\infty=0\,. 
$$ 
Returning to the original variables $(t,x)$, we obtain 
$$ 
\lim_{t\to\infty} \left\| (1+(q-1)t)^{1/(q-1)} u(t) - 
v_*\right\|_\infty=0\,, 
$$ 
whence 
\beqn 
\label{c8} 
\lim_{t\to\infty} \left\| t^{1/(q-1)} u(t) - (q-1)^{-1/(q-1)}\ 
v_*\right\|_\infty=0\,. 
\eeqn

It remains to identify $v_*$, or equivalently $V_*$. Since the latter 
is a viscosity solution to $|\nabla V_*|-1=0$ in $\mathcal{P}_\infty$ 
with $V_*=0$ on $\partial\mathcal{P}_\infty$ by the previous analysis, 
we infer from \cite[Remark~II.5.10]{BdCD97} that 
$$ 
V_*(x) = \inf_{y\in\partial\mathcal{P}_\infty}{\left\{ |y-x| \right\}} 
= \inf_{y\in\RR^N\setminus\mathcal{P}_\infty}{\left\{ |y-x| \right\}} 
\;\;\mbox{ for }\;\; x\in \mathcal{P}_\infty\,. 
$$ 
Consequently, $v_*=V_\infty$, the latter being defined by (\ref{a10}), 
and (\ref{a9}) follows from (\ref{c8}). \qed

\section{Invariant and moving positivity sets}\label{ips}

We first show that, if $u_0$ vanishes sufficiently rapidly near a 
point $x_0$, then $u(t,x_0)=0$ for all $t\ge 0$.

\begin{lem} 
\label{led1} 
Assume that there are $x_0\in\RR^N$ and $\delta>0$ such that 
\beqn 
\label{d2} 
u_0(x) \le a_0\ |x-x_0|^{(p-q)/(p-1-q)} \;\;\mbox{ for }\;\; x\in 
B(x_0,\delta)\,, 
\eeqn 
\beqn 
\label{d3} 
\| u_0\|_\infty \le a_0\ \delta^{(p-q)/(p-q-1)}\,, 
\eeqn 
the parameter $a_0$ being defined by {\rm (\ref{a16})}. Then, 
\beqn 
\label{d4} 
u(t,x_0)=0 \;\;\mbox{ for }\;\; t\ge 0\,. 
\eeqn 
\end{lem}

\noindent\textbf{Proof.} We adapt an argument from 
\cite[Theorem~8.2]{Kn77} and put $S_1(x) := a_0\ 
|x-x_0|^{(p-q)/(p-1-q)}$ for $x\in\RR^N$. An easy computation shows 
that 
$$ 
- \Delta_p S_1(x) + |\nabla S_1(x)|^q = 0 \;\;\mbox{ for }\;\; 
x\in\RR^N\,. 
$$ 
Furthermore, we have $u_0(x)\le S_1(x)$ for $x\in B(x_0,\delta)$ by 
(\ref{d2}) and $u(t,x)\le\|u_0\|_\infty \le S_1(x)$ for $(t,x)\in 
[0,\infty)\times\partial B(x_0,\delta)$ by (\ref{a8a}) and 
(\ref{d3}). The comparison principle \cite[Theorem~8.2]{CIL92} then 
entails that $u(t,x)\le S_1(x)$ for $(t,x)\in [0,\infty)\times 
B(x_0,\delta)$. Consequently, $u(t,x_0)\le S_1(x_0)=0$ for $t\ge 0$ 
and the lemma is proved.  \qed

\medskip

\noindent\textbf{Proof of Theorem~\ref{tha3}.} We first consider 
$x_0\in\partial\mathcal{P}_0$. Owing to (\ref{a14}) and (\ref{a15}) we 
are in a position to apply Lemma~\ref{led1} to deduce that 
$u(t,x_0)=0$ for all $t\ge 0$. We next consider $x_0$ 
lying outside the closure 
of $\mathcal{P}_0$ and $x\in B(x_0,\delta)$. Either $x\not\in 
\mathcal{P}_0$ and $u_0(x)=0\le a_0\ |x-x_0|^{(p-q/(p-1-q)}$. Or $x\in 
\mathcal{P}_0$ and there exists $\vartheta\in (0,1)$ such that 
$x_\vartheta:=\vartheta x + (1-\vartheta) x_0$ belongs to 
$\partial\mathcal{P}_0$. Then $|x-x_\vartheta|=(1-\vartheta)\ 
|x-x_0|\le\delta$ and we infer from (\ref{a14}) that $u_0(x)\le a_0\ 
|x-x_\vartheta|^{(p-q/(p-1-q)} \le a_0\ |x-x_0|^{(p-q/(p-1-q)}$. The 
assumptions of Lemma~\ref{led1} are then fulfilled by $u_0$ at $x_0$ 
and we conclude again that $u(t,x_0)=0$ for all $t\ge 0$.

We have thus shown that $\mathcal{P}(t)\subset\mathcal{P}_0$ for $t\ge 
0$. The opposite inclusion readily follows from the time monotonicity 
of the positivity set established in Proposition~\ref{prb2}. \qed

\medskip

We now turn to the proof of Theorem~\ref{tha4} and first establish the 
following result, which is in the spirit of 
\cite[Proposition~3.2]{AC83}.

\begin{lem}\label{led2} 
Consider $q\in (1,q_2)$, $q_2$ being defined in {\rm (\ref{a18})}. 
There are positive real numbers $C_4$ and $C_5$ depending only on 
$N$, $p$, and $q$ such that, if $u_0$ is an initial condition 
fulfilling {\rm (\ref{a5})} and 
$$ 
M:= \int_{B(0,1)} u_0(x)\ dx \ge C_5\,, 
$$ 
then the corresponding solution $u$ to {\rm (\ref{a1}), (\ref{a2})} 
satisfies 
$$ 
u(1,0) \ge C_4\ M^{p\eta} \;\;\mbox{ with }\;\; \eta:= 
\frac{1}{N(p-2)+p}\,. 
$$ 
\end{lem}

\noindent\textbf{Proof.} We adapt the proof given in \cite{Va06} for 
the porous medium equation, the main difference being that the 
$L^1$-norm of $u$ is not constant in our case. In the following, 
we denote by $C_i$, $i\ge 6$, positive constants depending only on 
$N$, $p$, and $q$.

Let us first assume that $\mbox{ supp }(u_0)\subset B(0,1)$, so 
that $\|u_0\|_1=M$. Denoting by $w$ the solution to the 
$p$-Laplacian equation $\partial_t w - \Delta_p w = 0$ in 
$Q$ with initial condition $w(0)=u_0$, the 
comparison principle entails that $u(t,x)\le w(t,x)$ for $(t,x)\in 
[0,\infty)\times\RR^N$. By \cite[Proposition~2.2]{KV88} and 
\cite[Proposition~1.3]{BtLxx} we have \bean \| u(t)\|_\infty & \le 
& \|v(t)\|_\infty \le C_6\ M^{p\eta}\ 
t^{-N\eta}\,,\\ 
u(t,x) &= & v(t,x) = 0 \;\;\mbox{ if }\;\; |x|\ge 1 + C_7\ 
M^{\eta(p-2)}\ t^\eta\,, \\ 
\|\nabla u(t)\|_\infty & \le & C\ M^{2\eta}\ t^{-\eta(N+1)}\,, 
\eean 
for $t>0$, the constants $C_6$ and $C_7$ depending only on $N$, $p$, 
and $q$. We next use the reflection argument of Aleksandrov as in the 
proof of \cite[Lemma~2.2]{AC83} to deduce from (\ref{a1}) that 
$$ 
u(t,0)\ge u(t,x) \;\;\mbox{ for }\;\; t\ge 0 \;\;\mbox{ and }\;\; 
|x|\ge 2\,. 
$$ 
We then infer from the above bounds that, if $C_7 M^{\eta(p-2)}>1$, 
\bean 
& & C_8\ u(1,0) \left[ \left( 1 + C_7\ M^{\eta(p-2)} \right)^N - 
2^N \right] \\ 
& \ge & \int_{\{|x|\ge 2\}} u(1,x)\ dx \\ 
& \ge & \|u(1)\|_1 - \int_{B(0,2)} u(1,x)\ dx \\ 
& \ge & M - \int_0^1 \int_{\RR^N} |\nabla u(s,x)|^q\ dxds - C_8\ 2^N\ 
\|u(1)\|_\infty \\ 
& \ge & M - C_8\ \int_0^1 \left( 1 + C_7\ M^{\eta(p-2)}\ s^\eta 
\right)^N\ \|\nabla u(s)\|_\infty^q\ ds - C_9\ M^{p\eta} \\ 
& \ge & M - C_{10}\ M^{N\eta(p-2)}\ M^{2\eta q)}\ \int_0^1 s^{-\eta q 
(N+1)}\ ds \\ 
& \ge & M - C_{11}\ \left( M^{\eta(N(p-2)+2q)} + M^{p\eta} \right)\,, 
\eean 
the assumption $q<p-(2N/(N+1))$ being used to obtain the last 
inequality. Consequently, if $C_7 M^{\eta(p-2)}>1$, 
$$ 
2^N\ C_8\ \left[ C_7^N\ M^{N\eta(p-2)} - 1 \right]\ u(1,0) \ge  M - 
C_{11}\ \left( M^{\eta(N(p-2)+2q)} + M^{p\eta} \right)\,. 
$$ 
Since $p\eta\in (0,1)$ and $\eta(N(p-2)+2q)\in (0,1)$ by 
(\ref{a18}) we readily conclude that there are $C_{12}$ and $C_{13}$ such 
that $u(1,0)\ge C_{12}\ M^{p\eta}$ provided $M\ge C_{13}$.

We next consider an arbitrary initial condition $u_0$ fulfilling 
(\ref{a5}). Then there is $\zeta\in\mathcal{C}_0^\infty(\RR^N)$ such 
that $0\le\zeta\le 1$, $\mbox{ supp }\zeta \subset B(0,1)$, and 
$$ 
\int_{B(0,1)} \zeta(x)\ u_0(x) = \frac{M}{2}\,. 
$$ 
Denoting by $\tilde{u}$ the solution to (\ref{a1}) with initial 
condition $\tilde{u}(0)=\zeta\ u_0$, it follows from the above 
analysis that $\tilde{u}(1,0)\ge C_{12}\ (M/2)^{p\eta}$ if $M\ge 2\ 
C_{13}$. As the comparison principle warrants that $u(1,0)\ge 
\tilde{u}(1,0)$, the expected result follows with $C_5=2\ C_{13}$ 
and $C_4=C_{12}\ 2^{-p\eta}$. \qed

\medskip

\noindent\textbf{Proof of Theorem~\ref{tha4}.} For $\lambda>0$, $t\in 
[0,\infty)$ and $x\in\RR^N$, we define 
$$ 
u_\lambda(t,x):= \lambda^{p-q}\ u\left( \lambda^{2q-p} t , 
x_0 + \lambda^{q-p+1} x \right)\,, 
$$ 
and observe that $u_\lambda$ also solves (\ref{a1}) with initial 
condition $u_\lambda(0)$. Furthermore, $u_\lambda(0)$ fulfils 
(\ref{a5}) and, if $\lambda\ge \delta^{-1/(p-1-q)}$, we infer from 
(\ref{a17}) that 
\bean 
\int_{B(0,1)} u_\lambda(0,x)\ dx & = & \lambda^{p-q+N(p-1-q)}\ 
\int_{B\left( x_0,\lambda^{q-p+1}\right)} u_0(x)\ dx \\ 
& \ge & A\ \lambda^{p-q+N(p-1-q)}\ 
\int_{B\left( x_0,\lambda^{q-p+1}\right)} |x-x_0|^{(p-q)/(p-1-q)}\ dx 
\\ 
& \ge & C_{14}\ A\,. 
\eean 
Therefore, if  $\lambda\ge \delta^{-1/(p-1-q)}$ and $A\ge A_q:= 
C_5/C_{14}$, we infer from Lemma~\ref{led2} that 
$$ 
\lambda^{p-q}\ u(\lambda^{2q-p},x_0) = u_\lambda(1,0) \ge C_4\ \left( 
\int_{B(0,1)} u_\lambda(0,x)\ dx \right)^{p\eta} > 0\,. 
$$ 
Consequently, $u(t,x_0)>0$ and thus $x_0\in\mathcal{P}(t)$ for $t\in \left( 
0,\delta^{(p-2q)/(p-1-q)} \right)$. We finally use the time 
monotonicity of $\mathcal{P}(t)$ established in Proposition~\ref{prb2} 
to complete the proof of Theorem~\ref{tha4}. \qed

\medskip

Combining Theorem~\ref{tha4} with a comparison argument allows us 
to extend Theorem~\ref{tha4} to $q\in [q_2,p-1)$ under stronger 
assumptions on the initial data.

\begin{prop}\label{prd3} 
Consider $x_0\in\partial\mathcal{P}_0$ and assume that there are 
$A>0$, $\delta>0$, and $r\in (1,q_2)$ such that \beqn \label{d8} 
u_0(x) \ge A\ |x-x_0|^{(p-r)/(p-1-r)} \;\;\mbox{ for }\;\; x\in 
B(x_0,\delta)\,. \eeqn Then, if $q\in [q_2,p-1)$, we have 
$x_0\in\mathcal{P}(t)$ for all $t>0$. 
\end{prop}

\noindent\textbf{Proof.} We put $\lambda_0:=1/\|\nabla 
u_0\|_\infty$ and $U(t,x) := \lambda_0^{p-q}\ u\left( 
\lambda_0^{2q-p} t, \lambda_0^{q-p+1} x \right)$ for $(t,x)\in 
[0,\infty)\times\RR^N$. Then $U$ solves (\ref{a1}) with initial 
condition $U_0:=U(0)$ and $\|\nabla U_0\|_\infty=1$. In addition, 
$X_0:=\lambda_0^{p-1-q} x_0$ clearly belongs to the boundary of 
the positivity set of $U_0$ and it follows from (\ref{d8}) that 
\beqn \label{d9} U_0(x) \le \lambda_0^{(q-1-r)/(p-1-r)}\ A\ 
|x-X_0|^{(p-r)/(p-1-r)} \;\;\mbox{ for }\;\; x\in B\left( 
X_0,\delta \lambda_0^{p-1-q} \right)\,. \eeqn Introducing 
$r_1:=(r+q_2)/2\in (r,q_2)$ we deduce from (\ref{d9}) that, if 
$\delta_1\in \left( 0,\delta \lambda_0^{p-1-q} \right)$ and $x\in 
B(X_0,\delta_1)$, \bean U_0(x) & \ge & 
\lambda_0^{(q-1-r)/(p-1-r)}\ A\ 
|x-X_0|^{(p-r_1)/(p-1-r_1)}\ |x-X_0|^{(r-r_1)/((p-1-r)(p-1-r_1))} \\ 
& \ge & \lambda_0^{(q-1-r)/(p-1-r)}\ 
\delta_1^{(r-r_1)/((p-1-r)(p-1-r_1))}\ A\ 
|x-X_0|^{(p-r_1)/(p-1-r_1)}\,, 
\eean 
whence 
\beqn 
\label{d10} 
U_0(x) \le A_1\ |x-X_0|^{(p-r_1)/(p-1-r_1)} \;\;\mbox{ for }\;\; x\in 
B\left( X_0,\delta_1 \right) 
\eeqn 
with $A_1 := \lambda_0^{(q-1-r)/(p-1-r)}\ 
\delta_1^{(r-r_1)/((p-1-r)(p-1-r_1))}\ A$. Furthermore, we can choose 
$\delta_1$ sufficiently small so that $A\ge A_{r_1}$, the constant 
$A_{r_1}$ being defined in Theorem~\ref{tha4}.

Now, on the one hand, it follows from (\ref{d10}) and 
Theorem~\ref{tha4} that the solution $\tilde{u}$ to 
\bean 
\partial_t\tilde{u} - \Delta_p\tilde{u} + 
\vert\nabla\tilde{u}\vert^{r_1} & = & 0\ \,, 
\quad (t,x)\in Q\,,\\ 
\tilde{u}(0) & = & U_0\,, \quad x\in\RR^N\,, 
\eean 
satisfies $\tilde{u}(t,X_0)>0$ for $t>0$. On the other hand, since 
$\|\nabla U_0\|_\infty=1$, we infer from (\ref{a1}) and (\ref{a8a}) 
that 
$$ 
\partial_t U - \Delta_p U + \vert\nabla U\vert^{r_1} \ge \partial_t U 
- \Delta_p U + \|\nabla U\|_\infty^{r_1-q}\ \vert\nabla U\vert^q \ge 
\partial_t U - \Delta_p U + \vert\nabla U\vert^q = 0 
$$ 
in $Q$. The comparison principle 
\cite[Theorem~8.2]{CIL92} then implies that $U(t,x)\ge \tilde{u}(t,x)$ 
for $(t,x)\in [0,\infty)\times\RR^N$. Consequently, 
$$ 
u(t,x_0) = \lambda_0^{q-p}\ U\left( \lambda_0^{p-2q} t,X_0 \right) \ge 
\lambda_0^{q-p}\ \tilde{u}\left( \lambda_0^{p-2q} t,X_0 \right) > 0 
$$ 
for $t>0$. \qed

\section{Comments, extensions, and open problems}

Theorem~\ref{tha2} gives a precise description of the behaviour 
for large times of solutions to (\ref{a1}) for $p>2$, $q\in 
(1,p-1)$, and non-negative and compactly supported initial data. 
It then gives rise to several related questions, concerning the 
situation when either the parameters $p$ and $q$ or the initial 
data do not fulfil the above conditions. 
 
 In particular, a natural question is whether a 
similar result is valid for non-negative initial data which decay 
rapidly at infinity but are not compactly supported. A first step 
in that direction would be to identify a class of initial data for 
which the solutions to (\ref{a1}), (\ref{a2}) satisfy the temporal 
decay estimates (\ref{a8b}).

Returning to compactly supported initial data, it is tempting to 
investigate what happens when $q$ reaches the boundary of the 
range $(1,p-1)$ analyzed here. On the one hand, the critical 
exponent $q=q_1=p-1$ with $p>2$ offers an interesting study of 
matched asymptotics that we plan to describe in a future 
publication. On the other hand, the critical case $q=1$ seems to 
be quite open and interesting, even in the semilinear case 
\cite{BRV96,BRV97}. 
 
For larger values of $q>p-1$ the situation is expected to be more 
classical: asymptotic simplification should take place for 
$q>q_*=p-(N/(N+1))$ in the sense that the large time behaviour 
will be governed by the diffusive part of the equation. For 
intermediate values of $q\in (p-1,q_*)$, the existence of very 
singular solutions has been established in \cite{Sh04a} and they 
are expected to describe the large time behaviour of compactly 
supported solutions. 
 
 Concerning waiting times, an interesting technical open problem 
is to figure out whether Theorem~\ref{tha4} is also true for $q\in 
[q_2,q_1)$ ($q_2$ is defined in formula (\ref{a18})).

As a final comment, let us mention that the techniques of this paper 
could possibly be applied to similar 
equations involving diffusion and absorption like 
\begin{equation*} 
\partial_t u= \Delta u^m- |\nabla u|^q \;\;\mbox{ in }\;\; Q\,, 
\end{equation*} 
but we have noticed that there are a number of difficulties, maybe 
technical. 
 
Equations with variable coefficients are also worth considering, 
as well as more general equations of the form \eqref{geneq}. In 
fact, the area of asymptotic behaviour of parabolic equations with 
variable coefficients is quite open even in the linear case. We 
refer to \cite{HPolSaf} for an interesting recent development for linear 
parabolic equations.

\

 \noindent {\sc Acknowledgment.} This work was started during a visit of 
PhL  to the Univ. Aut\'onoma de Madrid supported by Spanish Project 
MTM-2005-08760-C02-01. JLV is partially supported by this project and 
by ESF Programme ``Global and geometric aspects of nonlinear 
partial differential equations".

\ 
 
{\sl Keywords:}  nonlinear parabolic equations, $p$-Laplacian 
equation, asymptotic patterns, localization, Hamilton-Jacobi 
equations, viscosity solutions.


\begin{thebibliography}{99} 
 
 
 
\bibitem{ATU04} Daniele Andreucci, Anatoli~F.~Tedeev, and Maura Ughi, 
\textit{The Cauchy problem for degenerate parabolic equations with 
source and damping}, Ukrainian Math. Bull. \textbf{1} (2004), 
1--23. 
 
 
\bibitem{Ant96} Stanislav N. Antontsev, \textit{Quasilinear parabolic 
equations with 
non-isotropic nonlinearities: space and time localization}. Energy 
methods in continuum mechanics (Oviedo, 1994),  1--12, Kluwer 
Acad. Publ., Dordrecht, 1996. 
 
 
\bibitem{AntDiaz89} Stanislav N. Antontsev and Jes{\'u}s Ildefonso Diaz, 
\textit{ On space or time localization of solutions of nonlinear 
elliptic or parabolic equations via energy methods.} Recent 
advances in nonlinear elliptic and parabolic problems (Nancy, 
1988), 3--14, Pitman Res. Notes Math. Ser., 208, Longman Sci. 
Tech., Harlow, 1989. 
 
 
\bibitem{Ar85} Donald G.  Aronson,  \textit{The porous medium equation}. 
Nonlinear diffusion problems (Montecatini Terme, 1985),  1--46, 
Lecture Notes in Math., 1224, Springer, Berlin, 1986. 
 
 
\bibitem{AC83} Donald G.~Aronson and Luis A.~Caffarelli, 
\textit{The initial trace of a solution of the porous medium 
equation}, Trans. Amer. Math. Soc. \textbf{280} (1983), 351--366. 
 
 
 
\bibitem{BdCD97} 
Martino Bardi and Italo Capuzzo-Dolcetta, \textit{Optimal Control 
and Viscosity Solutions of Hamilton-Jacobi-Bellman Equations}, 
Systems Control Found. Appl., Birkh\"auser, Boston, 1997. 
 
 
 
\bibitem{Bl85} 
Guy Barles, \textit{Asymptotic behavior of viscosity solutions of 
first order Hamilton-Jacobi equations}, Ricerche Mat. \textbf{34} 
(1985), 227--260. 
 
 
 
\bibitem{Bl94} 
Guy Barles, \textit{Solutions de Viscosit\'e des Equations 
d'Hamilton-Jacobi}, Math\'ematiques \& Applications \textbf{17}, 
Springer-Verlag, Berlin, 1994. 
 
 
 
\bibitem{BlP88} 
Guy Barles and Beno\^\i t Perthame, \textit{Exit time problems in 
optimal control and vanishing viscosity method}, SIAM J. Control 
Optim. \textbf{26} (1988), 1133--1148. 
 
 
 
\bibitem{BtLxx} Jean-Philippe Bartier and Philippe Lauren\c cot, 
\textit{Gradient estimates for a degenerate parabolic equation 
with gradient absorption and applications}, (submitted). 
 
 
 
\bibitem{BKL04} 
Sa\"\i d Benachour, Grzegorz Karch, and Philippe Lauren\c cot, 
\textit{Asymptotic profiles of solutions to viscous 
Hamilton-Jacobi equations}, J. Math. Pures Appl. (9) \textbf{83} 
(2004), 1275--1308. 
 
 
 
\bibitem{BLSS02} 
Sa\"\i d Benachour, Philippe Lauren\c cot, Didier Schmitt, and 
Philippe Souplet, \textit{Extinction and non-extinction for 
viscous Hamilton-Jacobi equations in $\RR^N$}, Asymptot. Anal. 
\textbf{31} (2002), 229--246. 
 
 
 
\bibitem{BRV96} 
Sa\"\i d Benachour, Bernard Roynette, and Pierre Vallois, 
\textit{Solutions fondamentales de $u_t - \frac{1}{2}\ u_{xx} = 
\pm 
 \vert u_x\vert$}, 
Ast\'erisque \textbf{236} (1996), 41--71. 
 
 
 
\bibitem{BRV97} 
Sa\"\i d Benachour, Bernard Roynette, and Pierre Vallois, 
\textit{Asymptotic estimates of solutions of $u_t - \frac{1}{2}\ 
\Delta u = -\vert\nabla u\vert$ in $\RR_+\times\RR^d$, $d\ge 2$}, 
J. Funct. Anal. \textbf{144} (1997), 301--324. 
 
 
 
\bibitem{BNP} Michiel Bertsch, Tokumori Nanbu, and Lambertus A. Peletier, 
\textit{Decay of solutions of a degenerate nonlinear diffusion 
equation}, Nonlinear Anal.  \textbf{6} (1982), 539--554. 
 
 
 
\bibitem{BKP}  Michiel Bertsch,  Robert Kersner, and Lambertus A. Peletier, 
\textit{Sur le comportement de la fronti\`ere libre dans une 
\'equation en th\'eorie de la filtration}, C. R. Acad. 
Sci. Paris S\'er.~I \textbf{295} (1982), 63--66. 
 
 
 
\bibitem{BGK03} 
Piotr Biler, Mohammed Guedda, and Grzegorz Karch, 
\textit{Asymptotic properties of   solutions of the viscous 
Hamilton-Jacobi equation}, J. Evolution Equations \textbf{4} 
(2004), 75--97. 
 
 
 
\bibitem{ChV96} Manuela Chaves and Juan Luis V{\'a}zquez, 
\textit{Nonuniqueness in nonlinear 
heat propagation: a heat wave coming from infinity}, Differential 
Integral Equations  \textbf{9} (1996), 447--464. 
 
 
 
\bibitem{ChV99} Manuela Chaves and Juan Luis V{\'a}zquez, 
\textit{Free boundary layer formation in nonlinear heat 
propagation}, Comm. Partial Differential Equations \textbf{24} 
(1999), 1945--1965. 
 
 
 
\bibitem{CIL92} 
Michael G.~Crandall, Hitoshi Ishii, and Pierre-Louis Lions, 
\textit{User's guide to viscosity solutions of second order 
partial differential equations}, Bull. Amer. Math. Soc. (N.S.) 
\textbf{27} (1992), 1--67. 
 
 
 
\bibitem{DV85} 
Jes\' us Ildefonso Diaz and Laurent V\'eron, \textit{Local vanishing 
properties of solutions of elliptic and parabolic quasilinear 
equations}, Trans. Amer. Math. Soc. \textbf{290} (1985), 787--814. 
 
 
 
\bibitem{GL07} Thierry Gallay and Philippe Lauren\c cot, 
\textit{Asymptotic behavior for a viscous Hamilton-Jacobi equation 
with critical exponent}, Indiana Univ. Math. J. \textbf{56} 
(2007), 459--479. 
 
 
 
\bibitem{Gi05} 
Brian H.~Gilding, \textit{The Cauchy problem for $u_t=\Delta 
u+\vert \nabla u\vert^q$, large-time behaviour}, J. Math. Pures 
Appl. (9) \textbf{84} (2005), 753--785. 
 
\bibitem{HPolSaf} Juraj H\'uska, Peter Pol\'a\v cik, and Mikhail V. Safonov, 
\textit{Harnack inequalities, exponential separation, and 
perturbations of principal Floquet bundles for linear parabolic 
equations}, Preprint, 2007. 
 
\bibitem{GB02} 
Abdelilah Gmira and Benyouness Bettioui, \textit{On the 
selfsimilar solutions of a diffusion convection equation}, NoDEA 
Nonlinear Differential Equations Appl. \textbf{9} (2002), 
277--294. 
\bibitem{I87} Hitoshi Ishii, 
\textit{A simple, direct proof of uniqueness for solutions of the 
Hamilton-Jacobi equations of eikonal type}, Proc. Amer. Math. Soc. 
\textbf{100} (1987), 247--251. 
 
 
 
\bibitem{Ka87} 
Anatolii S.~Kalashnikov, \textit{Some problems of the qualitative 
theory of non-linear degenerate second-order parabolic equations}, 
Russian Math. Surveys \textbf{42} (1987), 169--222. 
 
 
 
\bibitem{KV88} 
Shoshana Kamin and Juan Luis V\'azquez, 
\textit{Fundamental solutions and asymptotic behaviour for the 
$p$-Laplacian equation}, 
Rev. Mat. Iberoamericana \textbf{4} (1988), 339--354. 
 
 
\bibitem{Kn77} 
Barry F.~Knerr \textit{The porous medium equation in one 
dimension}, Trans. Amer. Math. Soc. \textbf{234} (1977), 381--415. 
 
 
\bibitem{LSU} { Olga A. Ladyzhenskaya, Vsevolod A.  Solonnikov, and Nina N. 
Uraltseva.} \textit{Linear and Quasilinear Equations of Parabolic 
Type\/}, Amer. Math. Soc., Providence, R.I. 1968. 
 
 
 
\bibitem{NR99} 
Gawtum Namah and Jean-Michel Roquejoffre, \textit{Remarks on the 
long time behaviour of the solutions of Hamilton-Jacobi 
equations}, Comm. Partial Differential Equations \textbf{24} 
(1999), 883--893. 
 
 
 
\bibitem{Ro01} 
Jean-Michel Roquejoffre, \textit{Convergence to steady states or 
periodic solutions in a class of Hamilton-Jacobi equations}, J. 
Math. Pures Appl. (9) \textbf{80} (2001), 85--104. 
 
 
 
\bibitem{Sh04a} 
Shi Peihu, \textit{Self-similar singular solution of a 
$p$-Laplacian evolution equation with gradient absorption term}, 
J. Partial Differential Equations \textbf{17} (2004), 369--383. 
 
 
 
\bibitem{V91} Juan Luis V\' azquez, {\it Singular solutions and 
asymptotic behaviour of nonlinear parabolic equations}, International 
Conference on Differential Equations, Vol. 1, 2 (Barcelona, 1991), 
234--249, World Sci. Publ., River Edge, NJ, 1993. 
 
 
\bibitem{Va06} 
Juan Luis V\'azquez, \textit{Smoothing and decay estimates for 
nonlinear diffusion equations. Equations of porous medium type}, 
Oxford Lecture Ser. Math. Appl. \textbf{33}, Oxford University 
Press, Oxford, 2006. 
 
 
 
\bibitem{V07} Juan Luis V\' azquez, {\it The porous medium equation. 
Mathematical theory}, Oxford Mathematical Monographs. 
The Clarendon Press, Oxford University Press, Oxford, 2007. 
 
 
\bibitem{Y96} Yuan Hongjun, \textit{Localization condition for a nonlinear 
diffusion equation}, Chinese J. Contemp. Math. \textbf{17} (1996), 45--58. 
 
 
\end{thebibliography}
\end{document}